\documentclass[12pt]{amsart}
\usepackage{graphicx}
\usepackage{amssymb}
\usepackage{amscd}
\usepackage{hyperref}
\usepackage{amsmath}
\newtheorem{theorem}{Theorem}[section]
\newtheorem{corollary}[theorem]{Corollary}
\newtheorem{lemma}[theorem]{Lemma}

\newtheorem{definition}[theorem]{Definition}

\newtheorem{rem}[theorem]{Remark}

\begin{document}
\title{The Levi Problem On Strongly Pseudoconvex $G$-Bundles}
\author{Joe J Perez}
\thanks{MR Classification numbers: 32E40; 32W05; 43A30}
\maketitle
\begin{abstract}
Let $G$ be a unimodular Lie group, $X$ a compact manifold with boundary, and $M$ the total space of a principal bundle $G\to M\to X$ so that $M$ is also a strongly pseudoconvex complex manifold. In this work, we show that if $G$ acts by holomorphic transformations satisfying a local property, then the space of square-integrable holomorphic functions on $M$ is infinite $G$-dimensional. \end{abstract}
% 
%%%%%%%%%%%%%%
\section{Introduction}
%%%%%%%%%%%%%
%
\subsection{Basic definitions and statement of main theorem} Let $M$ be a complex manifold with nonempty smooth boundary $bM$,
  $\bar M=M\cup bM$, so that $M$ is the interior of $\bar M$, 
and ${\rm dim}_{\mathbb C}(M)=n$.  We will also assume for simplicity that $\bar M$ is a closed subset in $\widetilde{M}$, a
complex neighborhood of $\bar M$ so that the complex structure on $\widetilde{M}$ extends that of $M$, and every point of $\bar M$ is an
interior point of $\widetilde{M}$.

Let us choose a smooth function $\rho :\widetilde{M}\to \mathbb R$ so that
\[ M=\{z\mid \rho(z)<0\}, \ \ bM = \{z\mid \rho(z)=0\},\]
\noindent 
and for all $x\in bM$, we have $d\rho(x)\neq 0$. 
For any $x\in bM$ define the {\it holomorphic tangent plane} to the boundary at $x$ by
\[ T^{\mathbb C}_{x}(bM) = \{w\in \mathbb C^{n}\mid \sum_{k=1}^{n}\left.\frac{\partial\rho}{\partial z_{k}}\right|_x w_{k}=0\}.\]
For $x\in bM$, define the Levi form $L_{x}$ by
\[ L_{x}(w,\bar w) = 
  \sum_{j, k=1}^{n} \left.\frac{\partial^{2}\rho}{\partial z_{j}\partial \bar z_{k}}\right|_x  w_{j}\bar w_{k},\quad (w \in T^{\mathbb C}_{z}(bM)).\]
Then $M$ is said to be {\it strongly pseudoconvex} if for every $x\in bM$, the form $L_{x}$ is positive definite. Since $\rho$ is real-valued, the Taylor expansion at $x$ of $\rho$ is
\begin{equation}\label{taylor}\rho (z) = \rho (x) + 2 {\mathfrak Re}\ f(z,x) + L_x(z-x, \bar z - \bar x) + \mathcal O(|z-x|^3),\quad (z\in\mathbb C^n)\end{equation}
\noindent
with
\begin{equation}\label{LeviPoly}f(z,x) = \sum_{k=1}^n \left.\frac{\partial \rho}{\partial z_{k}}\right|_x (z_{k} - x_{k}) + \frac{1}{2} \sum_{jk=1}^n  \left.\frac{\partial ^2 \rho}{\partial z_{j} \partial z_{k}}\right|_x(z_{j}-x_{j})(z_{k}-x_{k}).
\end{equation}
This $f$ is holomorphic in $M\cap U_x$, with $U_x$ a small neighborhood of $x$, and vanishes only at $x$. We will see why later, but the point is the positivity of $L_x$. It happens that negative powers of $f$ are also holomorphic in the neighborhood and blow up at $x$. A question we will address in this paper is under which conditions we can correct this locally defined $f^\tau$ to obtain a global holomorphic function on $M$. In those conditions we will also say something about the size of $\mathcal O(M)$, the space of all holomorphic functions on $M$. 

The motivation behind this question is described in detail in \cite{Si}.  Early in the study of functions of several complex variables, the idea of pseudoconvexity arose in \cite{L}.

A point $x\in bM$ is called a {\it peak point for} $\mathcal O(M)$ if there exists an $f\in\mathcal O(M)$ unbounded on any neighborhood of $x$ and bounded in the complement of that neighborhood.
 
The Oka-Grauert theorem \cite{Gr} asserts that if $\bar M\subset \mathbb C^{n}$ is compact, has nonempty boundary, and is strongly pseudoconvex, then every point of the boundary is a peak point for $\mathcal O(M)$.  One way of proving this theorem and its variants will be described in this introduction.  
  
A point $x\in bM$ is called a {\it local peak point for} $\mathcal O(M)$ if there exist a function $f\in\mathcal O(M)$ and a neighborhood $V$ of $x$ in $M$ such that $f$ is unbounded on $V$, but bounded on $V\setminus U_x$ for any neighborhood $U_x$ of $x$ in $M$.  It was proven in \cite{GHS} that if $M$ is a strongly pseudoconvex complex manifold admitting a free cocompact holomorphic action of a discrete group, then every point in the boundary of $M$ is a local peak point for $L^2\cap\mathcal O(M)$, necessarily nontrivial. 
    
The goal of the present work is to extend this last  result (that $L^2\cap\mathcal O$ be nontrivial) from \cite{GHS} to general unimodular Lie group bundles. With a technical assumption (that we call amenability) on the local properties of convolutions of functions $f^\tau$ we will demonstrate
\begin{theorem} Assume that $G$ is a unimodular Lie group and $G\to M\to X$ a principal $G$-bundle.  Assume further that the total space $M$ is a strongly pseudoconvex complex manifold on which $G$ acts amenably by holomorphic transformations and that $X$ is compact. Then
\[\dim_G L^2\mathcal O(M)=\infty.\]
\end{theorem}
It is natural to assume the unimodularity of $G$ in this context and not only because it is an important tool in our formalism. In fact, \cite{GHS} contains a $G$-bundle with {\it non}unimodular structure group having $L^2\cap\mathcal O(M)=\{0\}$. However, unimodularity is not the whole story. The same paper also describes a $G$-manifold with nonunimodular structure group and many holomorphic functions.

Another word on the relationship between the results in \cite{GHS} and ours: if in our setting the structure group $G$ possesses a cocompact discrete subgroup $\Gamma$, then \cite{GHS} is applicable and one obtains local peak points and nontriviality solving a reduced problem:
\begin{equation}\label{reduce}G\to M\to X \quad\leadsto\quad \Gamma\to M\to (X\times (G/\Gamma)).\end{equation}
Generically, however, it is not the case that a unimodular Lie group have such a subgroup, {\it cf.} \cite{M}.
 
It does not seem to us that the methods in the present paper always allow direct construction of unbounded holomorphic functions. Still, when holomorphic functions can be constructed here, they are not smooth in the boundary. In particular, they will not possess holomorphic extensions beyond the boundary and so remain in the spirit of the early investigations of pseudoconvexity. A natural source of examples of complex manifolds satisfying the hypotheses (except perhaps amenability) are the Grauert tubes of unimodular Lie groups $G$ and of real-analytic manifolds of the form $K\times G$ with $K$ compact. 
\subsection{Compact case} We begin by reviewing the case when $M$ is compact, modifying the argument used in \cite{FK} to conform to our method.  We discuss the construction of holomorphic functions with peak points because it turns out that it is essentially our method of constructing any clearly nontrivial holomorphic functions in subsequent sections. 
Suppose $M$ is a compact complex manifold whose boundary is strongly pseudoconvex and for a point $x\in bM$, we want a holomorphic function blowing up at $x$.  Define the antiholomorphic exterior derivative $\bar\partial:\Lambda^{0,0}\to \Lambda^{0,1}$ in local coordinates $(z_k)_k$ by $\bar\partial u=\sum \frac{\partial u}{\partial \bar z_{k}}d\bar z_{k}$. If it can be established that
\[\bar\partial u=\phi\] 
has a smooth solution $u$ whenever $\phi$ is a smooth antiholomorphic one-form that satisfies the compatibility condition $\bar\partial\phi=0$, then we may construct the function desired.  The first step is to use the pseudoconvexity property of the boundary to construct a function $f$, holomorphic in a neighborhood $U_x$ of $x$, that blows up just at $x$, as indicated before.  Next, we can take a smooth function $\chi$ with support in $U_x$ that is identically equal 1 close to $x$.  Extending $\chi f$ by zero on the rest of $M$, we obtain a function, which we also call $\chi f$, defined everywhere and smooth away from $x$.  Furthermore, $\bar\partial(\chi f)=(\bar\partial\chi)f = 0$ near $x$, so $\bar\partial\chi f$ can be extended smoothly to the boundary.  If we can now find a smooth solution to $\bar\partial u = \bar\partial\chi f$, then $\chi f-u$ is holomorphic and must blow up at $x$ since $u$ is smooth up to the boundary.
Let us describe the construction of solutions $u\in L^{2}(M)$ to $\bar\partial u=\phi$ with $\phi \in L^{2}(M,\Lambda^{0,1})$, $\bar\partial \phi =0$.  Note that solutions will only be determined modulo the kernel of $\bar\partial$ consisting of square-integrable holomorphic functions.  Also, it is preferable to deal with self-adjoint operators, so since the Hilbert space adjoint $\bar \partial ^{*}$ of $\bar\partial$ satisfies $\overline{{\rm im}\bar\partial^*}=({\rm ker}\bar\partial)^{\perp}$,
it is sufficient to seek $u$ of the form $u=\bar\partial^{*}v$ satisfying
\begin{equation}\label{boh}\bar\partial\bar\partial^{*}v =\phi.
\end{equation}
This is a self-adjoint operator.  In order to eliminate the compatibility 
condition on $\phi$ (and obtain an operator related to the Dolbeault cohomology of $M$ also) let us add a term $\bar\partial^{*}\bar\partial v$, thus obtaining 
\begin{equation}\label{KL}
(\bar\partial\bar\partial^{*}+\bar\partial^{*}\bar\partial)v=\phi,
\end{equation}
\noindent
where $\phi$ need not be assumed to satisfy $\bar\partial \phi=0$.  Define the operator $\square =\bar\partial\bar\partial^{*}+\bar\partial^{*}\bar\partial$.  Notice that when $\bar\partial \phi =0$ is true, equation \eqref{KL} reduces to equation (\ref{boh}) because applying $\bar\partial$ to equation \eqref{KL} gives $\bar\partial\bar\partial^{*}\bar\partial v =0$ which in turn implies 
\[0=\langle\bar\partial\bar\partial^{*}\bar\partial v, 
\bar\partial v \rangle= \|\bar\partial^{*}\bar\partial v\|_{L^{2}(M)}^{2}.
\]
Thus the new term in equation \eqref{KL} vanishes when the compatibility condition holds. 
So it is enough to prove the solvability of the equation \eqref{KL}. But in fact, as we will see below, it suffices to prove that the operator $\square$ is {\it Fredholm}, {\it i.e.} the spaces ${\rm ker}\square$ and ${\rm coker}\square$ are finite-dimensional.
 
The equation $\square u=\phi$ is a noncoercive boundary value problem.  It has been shown \cite{K,FK} that on its domain in the antiholomorphic $q$-forms, when $q>0$, the operator $\square + 1$ has the following regularity property.  Let $\zeta, \zeta_{1}$ be smooth cutoff functions for which $\zeta_{1}=1$ on ${\rm supp}(\zeta)$ and let $H^{s}(M,\Lambda^{0,q})$ be the integer Sobolev space of sections in $\Lambda^{0,q}$ over $M$.  Then $\square v +v\in H^{s}_{\rm loc}(M,\Lambda^{0,q})$ implies $v\in H^{s+1}_{\rm loc}(M,\Lambda^{0,q})$ and there exist constants $C_{s}$ so that 
\begin{equation}\label{prima}\|\zeta v\|_{H^{s+1}(M)}\le 
  C_{s}\left(\|\zeta_{1}(\square + 1)v\|_{H^{s}(M)}+ \|(\square + 1)v\|_{L^{2}(M)}\right) \end{equation}
\noindent
uniformly in $v$.
These inequalities imply that the operator $(\square +1)^{-1}$ is bounded from $L^{2}(M,\Lambda^{0,q})$ to $H^{1}(M,\Lambda^{0,q})$ and so by Rellich's theorem is a compact operator in $L^{2}(M,\Lambda^{0,q})$ because $M$ is compact. Classical results of functional analysis allow one to conclude that $\square$ has discrete spectrum with no finite limit point and each eigenvalue has finite multiplicity.  Hence $\square$ has finite-dimensional kernel and cokernel and closed image ({\it i.e.} it is a Fredholm operator).   
  
Now, one can solve equation \eqref{KL} for all $\phi$ orthogonal to the finite-dimensional kernel.  As $\chi f$ is unbounded, raising $f$ to arbitrarily high powers generates linearly independent functions, still holomorphic in a neighborhood of $x$.  Further, since the $\chi f^{m}$ have compact support, $\bar\partial$ is injective on the vector space generated by $\{\chi f^{m}\mid m=1\dots N\}$. It follows that for $N$ sufficiently large, the image of this space under $\bar\partial$ intersects the image of $\square$ nontrivially:
\[ Q_{N}= {\rm im}\square\cap{\rm span}_{\mathbb C}\{\bar\partial \chi f^{m}\mid m=1\dots N\}\neq \{0\}.\]  
This, together with the fact that $Q_N\subset{\rm im}\bar\partial$ implies that $\bar\partial\bar\partial^{*} u= \phi$ can be solved for some $\phi\in Q_{N}$.  Since all the forms $\bar\partial \chi f^{m}$ are smooth, this $\phi$ will be smooth and so we proceed as indicated above.
\subsection{Regular coverings} As we have mentioned, in \cite{GHS} it was established that all boundary points are local peak points when $M$ is strongly pseudoconvex and admits a free cocompact action of a discrete group $\Gamma$ by holomorphic transformations. The proof above fails because when $M$ is not compact, Rellich's theorem no longer holds, so the dimension of the kernel and/or cokernel of $\square$ may be infinite-dimensional and the image of $\square$ may not be closed.  The {\it von Neumann dimension} of invariant subspaces of $L^{2}(\Gamma)$ is used in order to measure the kernel and cokernel of $\square$ in this setting as well as to measure the images of $\square$'s spectral projections.  We describe this briefly.  For a discrete group $\Gamma$, one forms 
\[ L^{2}(\Gamma)=\{\xi:\Gamma \to \mathbb C\mid 
\sum_{\gamma\in \Gamma} |\xi(\gamma)|^{2}<\infty\}.\]
This is a Hilbert space with inner product $\langle \xi,\eta\rangle_{L^{2}(\Gamma)} = 
\sum_{\gamma\in \Gamma} \xi(\gamma)\bar\eta(\gamma)$ and norm $\|\xi\|_{L^{2}(\Gamma)}^{2}=\langle\xi,\xi\rangle$.  Now, $\Gamma$ acts in 
$L^{2}(\Gamma)$ by right translations $R_{\gamma}$, $\gamma \in \Gamma$, defined by 
\[(R_{\gamma}\xi)(\alpha)= \xi(\alpha\gamma).\]
Clearly, $R_{\gamma}$ is a unitary operator.  A closed subspace $L\subset L^{2}(\Gamma)$ is called {\it invariant} if it is invariant with respect to $R_{\gamma}$ for all $\gamma \in \Gamma$.  It is true that if, in addition, our invariant subspace $L$ is closed, then $L$ is the image of a bounded left-convolution operator on the group:
\[L= {\rm im}L_{h} \quad {\rm where}\quad (L_{h}\xi)(\alpha)=
\sum_{\gamma \in \Gamma} h(\gamma) \xi(\gamma^{-1}\alpha)\]
\noindent
where $h:\Gamma\to\mathbb C$ is called a convolution kernel. 
Furthermore, one can choose $h$ so that $L_{h}$ is a self-adjoint projection: $L_{h}=L_{h}^{*}=L_{h}^{2}$.  Here the adjoint $L_{h}^{*}$ is defined by $\langle L_{h}^{*}\xi,\eta\rangle_{L^{2}(\Gamma)}=\langle\xi, L_{h}\eta\rangle_{L^{2}(\Gamma)}$ for all $\xi, \eta \in L^{2}(\Gamma)$.  

Defining $\mathcal B(L^{2}(\Gamma))$ to be the continuous linear operators in $L^{2}(\Gamma)$ and 
\[\mathcal L_{\Gamma}=\{L_{h}\mid h:\Gamma \to \mathbb C \ {\rm and}\ L_{h}\in \mathcal B(L^{2}(\Gamma))\}\]
\noindent
we see that $\mathcal L_{\Gamma}$ consists of all operators in $\mathcal B(L^{2}(\Gamma))$ commuting with the right-translations.  Von Neumann's bicommutant theorem then gives that $\mathcal L_{\Gamma}$ is a von Neumann algebra. On $\mathcal L_{\Gamma}$ there is a trace defined by 
\begin{equation}\label{gamtr}{\rm tr}_{\Gamma}(L_{h}) = h(e)\end{equation}
\noindent
and for a right-invariant subspace $L={\rm im}L_{h}$ with $L_{h}$ a self-adjoint projection, we define its $\Gamma${\it -dimension}
\[{\rm dim}_{\Gamma}(L)={\rm tr}_{\Gamma}(L_{h})=h(e).\]
Notice that since the identity in $\mathcal B (L^{2}(\Gamma))$ 
is convolution with $\delta$, the characteristic function of the identity, 
${\rm dim}_{\Gamma}(L^{2}(\Gamma))={\rm tr}_{\Gamma}(L_{\delta})=\delta(e)=1$, though of course ${\rm dim}_{\mathbb C}(L^{2}(\Gamma))=\infty$ for infinite groups.  

Next, when $\Gamma$ acts freely on a manifold $M$ with compact quotient, $X$, one decomposes the Hilbert space $L^{2}(M)\cong L^{2}(\Gamma)\otimes L^{2}(X)$ and defines a trace
\[{\rm Tr}_{\Gamma}={\rm tr}_{\Gamma}\otimes {\rm Tr}_{\mathcal B (L^{2}(X))}\] 
\noindent
on the invariant operators. It is with the corresponding dimension that closed, invariant subspaces of $L^{2}(M)$ are measured.  In \cite{GHS}, it is shown that a variant of Kohn's inequality \eqref{prima} implies that the kernel of $\square$ is finite-dimensional in this sense, though infinite-dimensional in the usual sense if nontrivial.  Moreover $\square$ is $\Gamma$-Fredholm in the sense that ${\rm im}\square$ contains a closed, $\Gamma$-invariant subspace of finite $\Gamma$-codimension. 

The operator $\square$ having the Fredholm property implies that the image of $\square$ intersected with
\begin{multline}\label{LN} L_{N}=\{\sum_{\gamma\in\Gamma}\sum_{m=1}^N c_{m,\gamma} (\bar\partial\chi f^m)(\cdot\ \gamma^{-1})\mid \sum_{m,\gamma}|c_{m,\gamma}|^2<\infty\} \\ 
\cong L^{2}(\Gamma)\otimes {\rm span}_{\mathbb C}\ \{\bar\partial\chi f, \bar\partial\chi f^{2}, \dots , \bar\partial\chi f^{N}\}\cong L^{2}(\Gamma)\otimes \mathbb C^{N}\end{multline}
\noindent
contains closed, invariant subspaces $Q$ of finite $\Gamma$-codimension in $L_{N}$.
Because $\bar\partial$ is injective on the span of the $\chi f^{m},\ m=1,2,\dots ,N$, we have that ${\rm dim}_{\Gamma}(L_{N})=N$.  As the kernel of $\square$ has finite $\Gamma$-dimension, the image of $\square$ contains closed, invariant subspaces of finite codimension, so the intersection ${\rm im}\square\cap L_{N}\subset L_{N}$ will be nontrivial if $N$ is sufficiently large.  Subsequently there exist closed, invariant nonempty subspaces $Q\subset {\rm im}\square\cap L_{N}$.  Picking a form $\phi\neq 0$ in this $Q$, one sees that it is smooth so $\square u=\phi$ is solvable and the rest of the argument is as previously described.

\subsection{$G$-bundles} In \cite{GHS} it is shown that $\square$ is $\Gamma$-Fredholm. In \cite{Per}, this theorem was adapted to the situation in which the discrete group $\Gamma$ is replaced by a unimodular Lie group $G$ (the reader is referred to \cite{Per} for the relevant definitions). 

For a unimodular group with its biinvariant measure fixed, the left- and right-convolutions $L_\Delta, R_\Delta$ by a distribution $\Delta$ on $G$ are defined as usual ({\it cf.} \S \ref{TGI}).  Also, the relevant von Neumann dimension is given by the trace ${\rm tr}_{G}$ on $\mathcal L_G \subset \mathcal B (L^2(G))$ agreeing with
\[{\rm tr}_G (L_{h}^*L_{h}) = \int_G |h(s)|^2 ds,\]    
whenever $L_{h}\in \mathcal B (L^{2}(G))$ and $h\in L^2(G).$  It is true that ${\rm tr}_G(A^* A)<\infty$ if and only if there is an $h \in L^2(G)$ for which $A= L_{h}\in \mathcal B(L^2(G))$.  

If we define $\tilde h(t) = \bar h(t^{-1})$, and  if $h_{j}, g_{j}\in L^2(G)$, $j=1,\dots, N$, then the operator $L_{k}= \sum_{1}^N L_{\tilde h_j} L_{g_j}$ is in ${\rm Dom}({\rm tr}_G)$. Furthermore, $k$ is continuous and ${\rm tr}_G(L_k) = k(e)$, agreeing with the discrete case Equation \eqref{gamtr}. We will outline the construction of the invariant trace ${\rm Tr}_G={\rm tr}_G \otimes {\rm Tr}_{\mathcal B(L^2(X))}$ below.

As we have suggested, in \cite{Per} we proved the following: {\it Assume that $G$ is a unimodular Lie group and $G\to M\to X$ a principal $G$-bundle.  Assume further that the total space $M$ is a strongly pseudoconvex complex manifold on which $G$ acts by holomorphic transformations and that $X$ is compact.  Then, for $q>0$, the operator $\square$ in $\Lambda^{p,q}(M)$ is $G$-Fredholm.} The $G$-Fredholm property is similar to the $\Gamma$-Fredholm property described above, {\it mutatis mutandis}.

In order to continue the program as in \cite{GHS}, the $\Gamma$-invariant spaces $L_N$ will need to be replaced by $G$-invariant versions. These spaces will be constructed similarly to the $L_N$ in Equation \eqref{LN}, namely by taking (some) convolutions of $\bar\partial\chi f$. As the bundle has a global right $G$ action, we may write convolutions on $M$ as we would on just $G$. The spaces replacing the $L_N$ will be $\{R_\Delta\bar\partial\chi f\mid P_\delta \Delta = \Delta\}\subset L^2(M)$ where $P_\delta$ is some projection in $\mathcal B(L^2(G))$ commuting with left-translations. With the $P_\delta$ chosen appropriately, these spaces are closed, smooth, right-invariant, and of arbitrarily large $G$-dimension, analogously to the $L_N$. Measuring these new spaces presents a new difficulty. In contrast to the discrete case in which $\dim_\Gamma L^2(\Gamma)=1$, there is now a complicated trace class for ${\rm tr}_G$. Our techniques here rely on methods developed in \cite{Per} with the exception of the new material in the present paper's Sections 2 and 3.  

Similarly to the previous cases, for $\delta>0$ sufficiently small, ${\rm im}\square\cap\{R_\Delta\bar\partial\chi f\mid P_\delta \Delta = \Delta\}$ will contain nontrivial closed subspaces so we proceed as usual except for one last contrast. 

In the compact case one constructs the function $\chi f-u$ and the singularity of $\chi f$ and smoothness of $u$ guarantee $\chi f-u\neq 0$. In the discrete group case, nothing changes in this respect. In the present situation we will be faced with the possibility that $R_\Delta\chi f$ be smooth to the boundary though $\chi f$ has a singularity there. In certain cases it is obvious that this cannot happen, but in others it is not (to us). We will handle a set of cases below in Section 4 and postpone a detailed discussion to a later paper. For now, let us say that our holomorphic action is {\it amenable} if there exists an $x\in bM$ so that if $f$ is the Levi polynomial at $x$, and $F$ is either some negative power of $f$ or the logarithm of $f$, then 1) $\chi f\in L^2(M)$, 2) $\|\chi F(\cdot,\xi)\|_{L^1(G)}^2<\infty$ for all $\xi\in X$, and 3) $R_\Delta\chi F\notin C^\infty(\bar M)$ for all $\Delta\in C^\infty(G)$ (we have chosen a local section $\xi:X\to M$ in the support of $\chi$). In the event that the action be amenable, we have our result arguing similarly as is done in the compact and covering space cases.

\subsection{Important examples} In \cite{GHS} a natural question is posed: is the cocompact unimodular group action relevant to the existence of holomorphic $L^2$-functions or is it just an toolmark of the method of proof? Now we might add another question and ask if the existence of holomorphic functions on a $G$-bundle has anything to do with amenability or if this is also just a useful tool in our proofs.

As we mentioned before, \cite{GHS} presents an example with the following properties. The complex dimension of $M$ is 2, $bM$ is strongly pseudoconvex, $G$ is a solvable nonunimodular connected Lie group, $\dim_{\mathbb R}G=3$, $G$ has a free action on $\bar M$ which is holomorphic on $M$, $\bar M/G=[-1,1]$, but $L^2\mathcal O(M)=\{0\}$.

The point here is that if we only impose bounded geometry conditions and uniformly strong pseudoconvexity, then the space of holomorphic $L^2$-functions may be trivial. 

Now a further property, the amenability, is involved that may or may not truly be relevant to the existence of holomorphic functions on such a manifold. Clearly more examples need be constructed and analyzed. We will postpone this for the future.

%%%%%%%%%%%%%%%

\subsection{Other approaches and results} Recent works on covering spaces extending \cite{GHS} are related to the Shafarevich conjecture (which asserts that the universal covering of a projective complex manifold is holomorphically convex) and can be found in \cite{Br1, Br2, Br3}. The paper \cite{TCM1} deals with the case in which $M$ is only assumed weakly pseudoconvex. Using cohomological techniques and holomorphic Morse inequalities, the authors obtain a lower bound for the $\Gamma$-dimension of the space of $L^2$ sections and upper bounds for the $\Gamma$-dimensions of the higher cohomology groups. In \cite{TCM2}, it is shown that the von Neumann dimension of the space of $L^2$ holomorphic sections is bounded below under weak curvature conditions on $M$.

In the present work, Section 2 contains a method for determining that a closed, $G$-invariant subspace of $L^2(M)$ be infinite $G$-dimensional. Section 3 describes a method of constructing large, smooth, invariant subspaces of $L^2(M)$. In Section 4 we construct local expressions for functions and convolutions and briefly discuss amenability.  In Section 5 we prove that $\dim_G L^2\mathcal O(M)=\infty$. Section 6 discusses a method by which the problem may be adjusted so as to give holomorphic functions with stronger singularities. 
\section{Paley-Wiener Theorems}
\noindent
This section is a small modification of a part of \cite{AL}.
\begin{definition}Let $M$ be a $G$-manifold with an invariant measure. For $f\in L^2(M)$, define $\langle f\rangle\subset L^2(M)$ to be the $L^2$-closure of the complex vector space generated by right-translates of $f$ by $G$. In symbols, 
\[\langle f\rangle = \overline{\left\{\sum_k^{\rm\small finite} \alpha_k f(\cdot\ t_k)\mid \alpha_k\in\mathbb C ,\ t_k\in G\right\}}^{\ L^2(M)}.\] 
\end{definition}  
%%%%%%%%%%%%%%%%%%%
\begin{theorem} \cite{AL} Let $G$ be a locally compact unimodular group containing a closed, noncompact, connected set. Let $f$ be in $L^2(G)$ such that $meas({\rm supp}(f))<meas(G)$ and such that there exists $h$ in $L^2(G)$ with $L_h f=f$. Then $f=0$, $m$-a.e.\end{theorem}
As we will need to recast a result from \cite{AL} in our language, we begin with an important fact about invariant operators in $L^2(G)$. As described in the introduction, on the von Neumann algebra $\mathcal L_G$ of bounded operators in $L^2(G)$ commuting with right-translations $R_t, t\in G$, there is a normal, faithful, semifinite trace ${\rm tr}_G$ agreeing with ${\rm tr}_G (L_{h}^*L_{h}) = \int_G |h(s)|^2 ds$, whenever this is defined ({\it cf.} \cite{T}). Using this invariant trace, we may define the dimensions $\dim_G$ of closed, right-invariant subspaces $L\subset L^2(G)$ as follows.  First, one notes that any such $L$ is the image of a self-adjoint projection $P_L$ in $\mathcal L_G$. As such, there exists a distribution $h$ on $G$ such that $P_L=L_h$. Then $\dim_G(L) = {\rm tr}_G(P_L) = \|h\|_{L^2(G)}^2$. 
\begin{corollary}\label{PW} Let $G\to M\to X$ be a principal $G$-bundle with $G$ a unimodular Lie group. If $0\neq h\in L^2(M)$ has sufficiently small support, then ${\rm dim}_G \langle h\rangle = \infty$.\end{corollary}
\begin{proof} Let the support of $h$ lie in a trivialization $G\times U$, $U\subset X$ of $M$ and choose a section so that we may write $h=h(t,x)$, $t\in G$, $x\in X$. Also let $P$ be a self-adjoint invariant projection whose image contains $\langle h\rangle$. By invariance 
\[PR_t h = R_t h\]
\noindent 
for any $t\in G$. By Lemma 1.2 of \cite{AL}, there exists a sequence $(t_k)_k\subset G$ for which the functions $(R_{t_k} h)_k$ are linearly independent and for which $S=\overline{\cup_k {\rm supp}(R_{t_k}h)}$ has finite measure. Denote by $\chi_S$ the characteristic function of $S$. The operator $u\mapsto \chi_S P u$ then has an infinite-dimensional eigenspace ${\rm span}\{R_{t_k}h\mid k\in\mathbb N\}$ corresponding to the eigenvalue one. We conclude that $\chi_S P$ must not be a compact operator.

\noindent
Let us compute the Hilbert-Schmidt norm of $\chi_S P$. Since $P$ is invariant, its representation in terms of its distributional kernel $\kappa$ takes the form  
\[(Pu)(t,x)=\int_{G\times X} ds dy\ \kappa(st^{-1};x,y)u(s,y).\]
If $(\psi_k)_k$ is an orthonormal basis for $L^2(X)$, the kernel of $\chi_S P$ can be expanded in a Fourier series
\[\chi_S(t)\kappa(st^{-1};x,y)=\chi_S(t)\sum_{kl} H_{kl}(st^{-1})\psi_k(x)\bar\psi_l(y).\]
\noindent
Since $(\psi_k\otimes\bar\psi_l)_{kl}$ forms an orthonormal basis for $L^{2}(X\times X)$, $H_{kl}$ is the $kl^{th}$ Fourier coefficient of $\kappa$ with respect to the decomposition $L^{2}(G\times X\times X)\cong \bigoplus_{kl}(L^{2}(G)\otimes\psi_k\otimes\bar\psi_l)$. We obtain
\[\|\chi_S P\|_{HS}^2 = \int_{G\times G}dsdt\ |\chi_S(t)|^2 \sum_{kl}|H_{kl}(st^{-1})|^2 \]\[= \sum_{kl}\|H_{kl}\|^2 \int_{G}dt\ |\chi_S(t)|^2 = meas(S) \sum_{kl}\|H_{kl}\|^2\]
\noindent
and conclude that $\sum_{kl}\|H_{kl}\|^2=+\infty$, for if not, we would have a Hilbert-Schmidt (and thus compact) operator $\chi_S P$ with an infinite-dimensional eigenspace corresponding to eigenvalue one.
We describe the invariant trace in $L^2(M)$, ({\it cf.} \cite{T}). Again using the orthonormal basis $(\psi_k)_k$ of $L^2(X)$, we have
\begin{equation}\label{decomp}
L^2(M) \cong L^{2}(G)\otimes L^2(X) \cong \bigoplus_{k\in \mathbb N} L^2(G)\otimes \psi_{k}.
\end{equation}
\noindent
Denoting by $P_k$ the projection onto the $k^{th}$ summand in \eqref{decomp}, we obtain a matrix representation of any operator $A\in \mathcal B(L^{2}(M))$ with elements $A_{kl} = P_k A P_l \in \mathcal B(L^{2}(G))$.  If $A\in \mathcal B(L^2(M))^G$, we recover the $H_{kl}$ from above as matrix elements 
\[ A \leftrightarrow [A_{kl}]_{kl}=[L_{H_{kl}}]_{kl}.\]
\noindent
The $G$-trace of such an operator is given by
\[{\rm Tr}_G(A) = \sum_k {\rm tr}_G(L_{H_{kk}}).\]
If $P$ is a self-adjoint projection, we compute ${\rm Tr}_{G}(P^{*}P)=\sum_{kl}{\rm tr}_{G}(L^{*}_{H_{kl}}L_{H_{kl}})=\sum_{kl}\|H_{kl}\|_{L^{2}(G)}^{2}$ by normality of ${\rm tr}_{G}$ and the definition of ${\rm tr}_G$. Thus $\dim_G \langle h\rangle = {\rm Tr}_G(P)=\sum_{kl}\|H_{kl}\|^2=\infty$.\end{proof}
% 
%%%%%%%%%%%%%%%%%%
\section{Smooth Invariant Closed Subspaces}
%%%%%%%%%%%%%%%%%%
%
\subsection{The group intrinsically}\label{TGI}
We gather some algebraic results. Define $\tilde\alpha(t)=\alpha(t^{-1})$ for any distributions $\alpha,\beta$ on $G$. The right-convolutions satisfy
\[(R_\alpha \beta)(t) \stackrel {\rm def} = \int_Gds\ \alpha(s)\beta(ts)=\int_Gds\ \beta(s)\alpha(t^{-1}s)=(R_\beta\alpha)(t^{-1}),\]
\noindent
so $R_\alpha \beta=\widetilde{R_\beta\alpha}$, and if $G$ is unimodular, then $\|R_\alpha \beta\|_{L^2(G)}= \|R_\beta \alpha\|_{L^2(G)}$. Using the definition $(L_s \alpha)(t)=\alpha(s^{-1}t)$, we obtain the identity 
\[(R_\alpha R_\beta \gamma)(t) = \int_Gds\ \alpha(s)\left[\int_G dr\ \beta(r) \gamma(tsr)\right]\]\[=\int_G dr\ \left[\int_G ds\ \alpha(s)\beta(s^{-1}r)\right]\gamma(tr) = (R_{[L_\alpha \beta]}\gamma)(t).\]
In this subsection, assume $H\in C^\infty_c(G)$ and consider $\langle H\rangle\subset L^2(G)$.  Any $g\in\langle H\rangle$ satisfies $g=\lim_m g_m$ with $g_m=R_{\Delta_m}H$ for some sequence $(\Delta_m)_m\subset C^\infty_c(G)$.
\noindent
Equivalently, $(g_m)_m$ is Cauchy, thus
\begin{equation}\label{cauchy}\|g_m -g_n\|=\|(R_{\Delta_m}-R_{\Delta_n})H\|=\|R_H(\Delta_m -\Delta_n)\|\longrightarrow 0.\end{equation}
%%%%%%%%%%%%%
\begin{definition}Let $R_H=U|R_H|$ be the polar decomposition of $R_H$, with $U$ a partial isometry, and let $|R_H|=\int_0^C\lambda dE_\lambda$ be the spectral decomposition of $|R_H|$. For $\delta\in [0,C]\cup \{0^+\}$, let $P_\delta=\int_\delta^C dE_\lambda$ and define 
\[\langle H\rangle_\delta = \{g\in\langle H\rangle\mid P_\delta U^*\tilde g=U^* \tilde g\}.\]
\end{definition}
\begin{lemma}\label{A}If $\delta>0$, then $g\in\langle H\rangle_\delta$ implies that $g=R_\Delta H$ for some $\Delta\in L^2(G)$. Consequently, $\langle H\rangle_\delta\subset H^\infty(G)$.\end{lemma}
\begin{proof}As in \eqref{cauchy}, let $R_{\Delta_m}H\to g\in\langle H\rangle_\delta$. Then $R_H \Delta_m\to \tilde g$ and
\[ U P_\delta U^*R_H \Delta_m\to  U P_\delta U^*\tilde g=\tilde g.\]
\noindent
The composition $P_\delta U^*R_H=P_\delta |R_H|=P_\delta |R_H|P_\delta$, when restricted to the orthogonal complement of ${\rm ker}P_\delta$, is an injection with bounded inverse, as is $U P_\delta U^*R_H$. Therefore there exists a Cauchy sequence $(\Delta_m')_m$ in $L^2(G)\ominus{\rm ker}P_\delta$ with limit $\Delta^g\in L^2(G)\ominus{\rm ker}P_\delta$ so that 
\[g=R_{\Delta^g} H.\]
\noindent
Noting that $\Delta^g\in L^2(G)$ for all $g\in\langle H\rangle_\delta$ and $H\in H_c^\infty(G)$, we have $\langle H\rangle_\delta\subset H^\infty(G)$. \end{proof}
\begin{rem}{\rm Since ${\rm im}|R_H|={\rm im}(R_H^*R_H)\subset C^\infty(G)$, we have ${\rm im}P_\delta\subset C^\infty(G)$ for all $\delta\in(0,C]$. Lemma \ref{A} and Corollary 6.4 of \cite{Per} provide that ${\rm dim}_G\langle H\rangle_\delta<\infty$ for $\delta>0$. The previous lemma gives that, if $\delta>0$, then $\langle H\rangle_\delta \subset \{R_\Delta H\mid \Delta\in {\rm im}P_\delta\}$. In fact, the spaces are equal:
}\end{rem}
\begin{lemma}\label{B} Let $|R_H|=\int_0^C\lambda dE_\lambda$ and $P_\delta=\int_\delta^C dE_\lambda$ as before. Then, for any $\delta>0$, we have $\langle H\rangle_\delta=\{R_\Delta H\mid \Delta\in{\rm im}P_\delta\}$.\end{lemma}
\begin{proof}For $\delta>0$, all $g\in\langle H\rangle_\delta$ satisfy
\[\tilde g=UP_\delta U^*\tilde g=UP_\delta U^*R_H\Delta^g=UP_\delta|R_H|\Delta^g=U|R_H|P_\delta\Delta^g=R_H P_\delta\Delta^g,\]
\noindent
so each $g\in\langle H\rangle_\delta$ is of the form $R_{\Delta^g}H$ for $\Delta^g\in{\rm im}P_\delta$. Conversely, if $\tilde g=R_H P_\delta\Delta^g$ for $\Delta^g\in{\rm im}P_\delta$, the above chain of equalities can be read right to left, obtaining $\tilde g=R_H P_\delta\Delta^g=UP_\delta U^*\tilde g$.\end{proof}
\begin{theorem}\label{C}For $\delta\in (0,C]$, the spaces $\langle H\rangle_\delta\subset\langle H\rangle$ are closed, smooth, right-invariant, and $\dim_G \langle H\rangle_\delta\to\infty$ as $\delta\to 0^+$.\end{theorem}
\begin{proof}The invariance condition on $\langle H\rangle_\delta$ is equivalent to the statement $g=R_\Delta H$ for $\Delta\in {\rm im}P_\delta$ if and only if
\[R_t g=R_t R_\Delta H = R_{[L_t \Delta]}H\in\langle H\rangle_\delta\qquad (t\in G).\]
\noindent
Since $P_\delta$ is a function of $R_H$, it commutes with all left-translations so $L_t \Delta\in {\rm im}P_\delta$. 
For a moment consider the case in which $\delta=0^+$. The projection $P_{0^+}$ is onto the closure of the image of $|R_H|=U^* R_H$ and so $P_{0^+}U^* = U^*$. The condition restricting $\langle H\rangle_\delta$, $P_{0^+} U^*\tilde g=U^* \tilde g$, is therefore vacuous, so $\langle H\rangle_{0^+}=\langle H\rangle$. By Corollary \ref{PW}, $\dim_G (\langle H\rangle_{0^+})=\infty$. Now, under the map $g\mapsto\tilde g$ we obtain an isomorphism 
\begin{equation}\label{antiisom}\langle H\rangle_\delta=\{R_{P_\delta\Delta}H\mid\Delta\in L^2(G)\}\cong\{R_H P_\delta\Delta\mid\Delta\in L^2(G)\}=\widetilde{\langle H\rangle}_\delta.\end{equation}
\noindent
Note that this isomorphism interchanges a left invariant subspace with a right invariant one. Let us see why they have the same $G$-dimension.  If $\langle H\rangle_\delta$ is the image of $L_h$, a self-adjoint projection, then $\widetilde{\langle H\rangle}_\delta$ is the image of the self-adjoint projection $R_h$. This is because $L_h g=g$ if and only if $R_h \tilde g =\tilde g$. We conclude that ${\rm tr}_G(L_h)={\rm tr}_G(R_h)=\|h\|_{L^2(G)}^2$, which implies that $\langle H\rangle_\delta$ and $\widetilde{\langle H\rangle}_\delta$ have the same $G$-dimension, though one is a left module and the other is a right module. 
In the polar decomposition $R_H=U|R_H|$, the partial isometry $U$ commutes with left-translations so, as left $G$-modules,
\[\widetilde{\langle H\rangle}_\delta=\{R_H P_\delta\Delta\mid\Delta\in L^2(G)\}\cong\{|R_H| P_\delta\Delta\mid\Delta\in L^2(G)\}\quad (\delta>0).\]
It is obvious that the right hand side of the above expression is equal to ${\rm im}P_\delta$, so combining these observations with Equation \eqref{antiisom}, we obtain ${\rm dim}_G \langle H\rangle_\delta = {\rm tr}_G(P_\delta)$, valid for $\delta>0$. 
Since $(P_\delta)_\delta$ are a spectral family,  $P_\delta\to P_{0^+}$ strongly (implying ${\rm dim}_G \langle H\rangle_\delta={\rm tr}_G P_\delta\to{\rm tr}_GP_{0^+}$), and so normality of the trace gives the result as soon as we obtain ${\rm tr}_G P_{0^+} = \infty$. To wit, $U\ {\rm im}P_{0^+}=\overline{{\rm im}U P_{0^+}}=\overline{{\rm im}R_H}$ and this has infinite $G$-dimension since it contains $H$. \end{proof}
\subsection{Actions}\label{actions} For a function $h\in C^\infty_c(\bar M)$ with small enough support, we may choose a section and write $h$ as a smooth function of $(t,x)\in G\times U$ where $U\subset X$. 
Since $M$ has a global right $G$-action, we may abbreviate a convolution by $\Delta$, $R_\Delta\otimes {\bf 1}_{L^2(X)}$, simply writing $R_\Delta$. We obtain an expression for $\|R_\Delta h\|_{L^2(M)}$ by first decomposing $h$ as in Equation \eqref{decomp}. With $H_k(t) = \langle h(t,\cdot),\psi_k\rangle_{L^2(X)}$, the function $h=\sum_k H_k\otimes\psi_k$ and 
\begin{equation}\label{dec}R_\Delta h=\sum_k (R_\Delta H_k)\otimes\psi_k \quad {\rm so}\quad  \|R_\Delta h\|_{L^2(M)}^2=\sum_k \|R_\Delta H_k\|_{L^2(G)}^2.\end{equation}
\begin{rem}\label{dls}{\rm Let $\delta>0$ and consider the decompositions $R_{H_k} = U_k |R_{H_k } |$, $|R_{H_k}|=\int_\delta^C \lambda dE^k_\lambda$ and the projections $P^k_\delta=\int_\delta^C dE^k_\lambda$. 
Then, for each $l\in\mathbb N$ for which $R_\Delta H_l\neq 0$ we have
\[\|R_\Delta h\|_{L^2(M)}^2=\sum_k \|R_\Delta H_k\|_{L^2(G)}^2\ge  \|R_\Delta H_l\|_{L^2(G)}^2 \ge \delta^2 \|\Delta\|_{L^2(G)}^2 \qquad (\Delta\in {\rm im}P^l_\delta).\]
\noindent
This implies that ${\rm im}P^l_\delta\ni\Delta\mapsto R_\Delta h$ is boundedly invertible as long as $R_\Delta H_l\neq 0$. Let us then take $\mathfrak D_\delta^l = {\rm im}P_\delta^l$ for $R_\Delta H_l\neq 0$ and define 
\[\langle h\rangle_{\delta,l}=\{R_\Delta h\mid \Delta\in\mathfrak D_\delta^l\}.\]}\end{rem}
\begin{lemma}\label{smoothbig}For $\delta>0$, the spaces $\langle h\rangle_{\delta,l}$ are closed, invariant, and smooth. Furthermore, $\dim_G \langle h\rangle_\delta<\infty$.\end{lemma}
\begin{proof}The previous remark and Lemma \ref{B} give that the space $\langle h\rangle_{\delta,l}$ is closed. For $\delta>0$, Lemma \ref{B} also provides that $\mathfrak D^l_\delta\subset C^\infty\cap L^2(G)$. Consider the estimate
\[\|R_\Delta h\|_{L^2(M)}^2 = \int_X dx\int_G dt\ \left|\int_Gds\ \Delta(s)h(ts,x) \right|^2 \le\|\Delta\|_{L^2(G)}^2\int_X dx\ \|h(\cdot,x) \|_{L^1(G)}^2  \]
\begin{equation}\label{L1}\lesssim \|\Delta\|_{L^2(G)}^2\left|\int_X dx\ \|h(\cdot,x)\|_{L^1(G)}\right|^2 = \|\Delta\|_{L^2(G)}^2\|h\|_{L^1(M)}^2,\end{equation}
\noindent
where $A\lesssim B$ means that for some $C>0$, $|A|\le C |B|$ uniformly. Recall the right-invariant Sobolev spaces as in \cite{Per}. There, the derivatives defining the spaces essentially commute with right translations. Thus the above estimate implies $\|R_\Delta h\|_{H^s(M)}^2 \lesssim \|\Delta\|_{L^2(G)}^2\|h\|_{L^1_s(M)}^2$ where $\|\cdot\|_{L^1_s(M)}$ is the $L^1$ norm of the derivatives up to order $s$. Since all derivatives of $h$ are in $L^1(G)$, we have $\langle h\rangle_{\delta,l}\subset H^\infty(\bar M)$. Corollary 6.4 of \cite{Per} states that if a space is closed, invariant, and in $H^\infty(\bar M)$, then it has finite $G$-dimension. 
\end{proof}
\begin{lemma}\label{bigspaces} As $\delta\to 0^+$, $\dim_G(\langle h\rangle_{\delta,l})\longrightarrow+\infty$.\end{lemma}
\begin{proof}By Equation \eqref{dec}, the space $\langle h\rangle_{\delta,l}$ has an orthogonal decomposition 
\[\langle h\rangle_{\delta,l} = \bigoplus_k \{R_\Delta H_k\mid \Delta\in\mathfrak D^l_\delta\}\otimes\psi_k .\]
\noindent
Now, since $H_l(t)=\langle h(t,\cdot),\psi_l\rangle_{L^2(X)}\in C^\infty_c(G)$, Theorem \ref{C} holds and provides that
\[\langle H_l\rangle_\delta\otimes\psi_l=\{R_\Delta H_l\mid \Delta\in\mathfrak D^l_\delta\}\otimes\psi_l\]
is a closed, invariant subspace of $\langle h\rangle_{\delta,l}$ whose $G$-dimension is unbounded as $\delta\to 0^+$.\end{proof}
%
%%%%%%%%%%%%%%%%%%%%%%%%%
\section{Levi's Function and Its Convolutions}
%%%%%%%%%%%%%%%%%%%%%%%%%
As discussed in the introduction, we need to know when convolutions of the singular functions gotten by taking the Levi polynomial to negative powers are not smooth in the boundary. We will not fully answer this question here but provide some tools and some simple examples, postponing a full analysis of the situation. We start with an analysis of $f$ itself. This first bit is in \cite{GHS}.

Without loss of generality (replacing $\rho$ by $e^{\lambda\rho}-1$ with sufficiently large $\lambda>0$) we may choose a defining function of $M$ so that the Levi form $L_x(w,\bar w)$ is positive for all nonzero $w\in\mathbb C^n$ (and not only for $w\in T^c_x(bM)$) and at all points $x\in bM$. Let us also assume that the defining function $\rho$ is constant on the orbits of points of $M$ and reconsider the Levi polynomial in Equation \eqref{LeviPoly}. The complex quadric hypersurface $S_x=\{z\mid f(z,x)=0\}$ has $T_x^{\mathbb C}(bM)$ as its tangent plane at $x$. The strong pseudoconvexity property implies that $\rho(z)>0$ if $f(z,x)=0$ and $z\neq x$ is close to $x$. This means that near $x$ the intersection of $S_x$ with $bM$ contains only $x$. The function $1/f(\cdot,x)$ is therefore holomorphic in $U\cap M$ (where $U$ is a neighborhood of $x$ in $\widetilde{M}$) and $x$ is its peak point. Since $\rho<0$ in $M$, \eqref{taylor} implies that $\mathfrak{Re}\  f(z,x)<0$ if $x\in bM$ and $z\in M$ is sufficiently close to $x$. It follows that we can choose a branch of $\log f(z,x)$
so that $g_x(z)=\log f(z,x)$ is a holomorphic function in
$z\in M\cap U_x$ where $U_x$ is a sufficiently small neighborhood of $x$ in $bM$. Consequently all powers of $f$ are also well-defined and holomorphic in a neighborhood of zero. Thus define
\[a = \left.\frac{\partial\rho}{\partial z_k}\right|_x, \qquad M = \frac{1}{2}\left.\frac{\partial^2\rho}{\partial z_k \partial z_l}\right|_x, \quad\quad f_{\tau}(z) = \left(\frac{1}{a\cdot z +
M z\cdot z}\right)^\tau\quad (\tau>0),\]
\noindent
and $f_{0}(z) = \log(a\cdot z+M z\cdot z)$ where $a\cdot b = \sum a_k b_k$. The functions $f_\tau$ are holomorphic in a neighborhood of $0$ and blow up only at $0$. 
\begin{lemma}\label{thanksMAS} Take coordinates in which $x$ is zero in the above. Then for $z$ sufficiently near zero in $\bar M$ there are constants $C,D>0$ so that
\[ C|z|^{2}\le |a\cdot z + M z\cdot z|\le D|z|.\]\end{lemma}
\begin{proof}This is true because 
\[ 2|a\cdot z + M z\cdot z|\ge -2\mathfrak{Re}\  (a\cdot z + M z\cdot z)\]
\[\ge \rho(z)-2\mathfrak{Re}\  (a\cdot z + M z\cdot z)=L_0(z,\bar z)+\mathcal O(|z|^3)\]
\noindent
and the Levi form has a smallest eigenvalue $\lambda>0$, so $L_0(z,\bar z)>\lambda|z|^2$. The other estimate is obvious.\end{proof}
Let $U$ be a neighborhood of a point $x=0$ of the boundary and choose a cut-off function $\chi\in C_c^\infty(U)$, so that $\chi=1$ in a neighborhood of $0$.  Locally defined functions cut off by $\chi$ will be considered extended by zero tacitly.
\begin{lemma}\label{l2} Let $\chi\in C^\infty_c(\bar M)$ with small support near zero. Then $\chi f_{\tau} \in L^p(M)$ whenever $\tau\in[0,n/p)$.\end{lemma}
\begin{proof}
For $\chi f_{\tau} \in L^p$ we only need $f_{\tau}\in L^p_{loc}$.  By Lemma \ref{thanksMAS}, there is a constant $C'$ so that, with $r=|z|$
\begin{equation}\int_{B_{\epsilon}} |f_{\tau}|^p dV \le C' \int_0^\epsilon
\frac{r^{2n-1}}{r^{2p\tau}}dr <\infty\end{equation}
\noindent
whenever $\tau<n/p$. The case of the logarithm is similar.\end{proof}
\begin{rem}\label{slice}{\rm Note that in the estimate \eqref{L1}, the quantity $ \|h(\cdot,x) \|_{L^1(G)}^2$ plays an important role. Later we will need a similar quantity $ \|\chi f_\tau(\cdot,x)\|_{L^1(G)}^2$.  The estimate above gives that for all $x\in X$, $\|\chi f_\tau(\cdot,x)\|_{L^1(G)}^2<\infty$ as long as $2\tau<\dim_{\mathbb R}G$}.\end{rem}
\begin{corollary}\label{L2} Let $d=\dim_{\mathbb R}G$, $\tau\in[0,d/2)$ and $\Delta\in L^2(G)$. Then $R_\Delta\chi f_\tau\in L^2(M)$.\end{corollary}
\begin{definition} Let $\xi:X\to M$ be a piecewise continuous section of $G\to M\stackrel p \to X$ so that $\xi|_{p({\rm supp}\chi)}$ is continuous. The action of $G$ on $M$ is called {\rm amenable} if there exists an $x\in bM$ and $\tau\ge0$ so that if $f$ is the Levi polynomial at $x$, then 1) $\chi f_\tau\in L^2(M)$, 2) $\|\chi f_\tau(\cdot,\xi)\|_{L^1(G)}^2<\infty$ for all $\xi\in X$, and 3) $R_\Delta\chi f_\tau\notin C^\infty(\bar M)$ for all nonzero $\Delta\in C^\infty(G)$.\end{definition}

%\begin{definition} Let $\tau\ge0$. We say that the $G$-action on $M$ is {\rm amenable} if there exists an $x\in bM$ so that if $f=f(\cdot,x)$ is the Levi polynomial with base point $x$, then $R_\Delta\chi f_\tau \notin C^\infty (\bar M)$ for any $\Delta\in C^\infty(G)$.\end{definition}
%
\subsection{Examples}
Here we will present some simple examples regarding amenability and the nonunimodular examples from \cite{GHS}. 
\subsubsection{The strip} Consider the tube of $\mathbb R$:
\[S=\{z\in\mathbb C\mid \mathfrak{Im}(z)\in(-1,1)\}.\]
The defining function is $\rho(z) = \mathfrak{Im}(z)^2-1$, and with the basepoint $x=a+ib\in bS$ for the Taylor expansion of $\rho$, the Levi polynomial is
\[f(z,\tilde\xi) = -i2b (z-\tilde\xi) - (z-\tilde\xi)^2.\]
Let us take the basepoint $x = i\leftrightarrow (a,b)=(0,1)$. Then, 
\[f(z) = f(z,i) = -2i(z-i)-(z-i)^2 = -(z^2+1).\]
The strip $S$ has an obvious translation invariance, 
\[\mathbb R\ni t:z\longmapsto z+t\]
and so convolutions of functions make sense on $S$. 
We can verify that $\mathfrak{Re}f<0$ close by the basepoint $i\in bS$, so $f$ can be raised to arbitrary real powers. Consider then the convolution
\[(\Delta\ast\chi f^\tau)(z) = \int_{\mathbb R}dt\ \Delta(t)\chi(z+t)f^\tau(z+t)\]
and its limits as $S\ni z\to bS$, or more specifically, $z\to i$. These are approximately
\[\lim_{z\to i}\int_{-\epsilon}^\epsilon dt\ [(z+t)^2+1]^\tau\]
\noindent
where $\epsilon\approx{\rm diam}({\rm supp}\chi)$ and we have ignored the details of $\Delta\in C^\infty(\mathbb R)$. These functions are well-behaved at the boundary, ({\it cf.} \cite{CS}, Theorem 2.1.3) and so this is not an amenable action. 

Of course, if we wanted to investigate $L^2\mathcal O$ in the strip, with the action of $\Gamma=\mathbb Z$ by translations along $\mathbb R$, we could use \cite{GHS}, which applies. In other words, the structure group $\mathbb R$ can be reduced to $\mathbb Z$. Also, in this case $L^2$ holomorphic functions can be obtained by explicit construction, taking {\it e.g.} $f(z)=1/(a^2+z^2)$ where $a>1$, or $\exp(-z^2)\log(z-i)$.

\subsubsection{Countably many strips} It may happen that $G$ is nonunimodular but there are many $L^2$-holomorphic functions on $M$. Consider the group $G\cong\mathbb Z\ltimes \mathbb R$ of matrices
\[ g=\left[
\begin{array}{ccc}
  2^n & x \\
  0 & 1       
\end{array}
\right]\qquad (n\in\mathbb Z,\ x\in\mathbb R)\]
\noindent
and let $M$ be the tube of $G$ consisting of all matrices

\[h=\left[
\begin{array}{ccc}
2^m & z  \\
0 & 1      
\end{array}
\right]
\qquad (m\in\mathbb Z,\ z\in\mathbb C,\ |\mathfrak{Im} z|<1).\]
Then $M$ is a disjoint countable union of strips
$\{z\in\mathbb C\mid |\mathfrak{Im} z|<1\}$. Consider $M$ as a (non connected)
complex manifold with boundary with $\bar M$ obtained by taking closure of each strip, so that $\bar M$ consists of the matrices of the same form with $|\mathfrak {Im}z|\le 1$. It is clear that $M$ is strongly pseudoconvex.

The action of $G$ on $M$ is obtained by left multiplication of the matrices: $g\cdot h=gh$. It amounts to interchanging the strips and real translations along the strips depending on the strip. This action is obviously holomorphic and free. It is easy to see that $\bar M/G=[-1,1]$.

Introducing the standard Euclidean metric on every strip (the metric induced by the standard
metric on $\mathbb C$), we obtain an invariant metric. The $L^2$ holomorphic functions on each strip constructed before, can be extended to $M$ by $0$. Here, questions of amenability are similar to the single strip case. 

%%%%%%%%%%%%%%%%%%%%%

\subsubsection{Trivial $L^2\mathcal O$} Let $\Omega\subset\mathbb C^2$ be the Siegel domain of the second kind, $\Omega = \{(z_1,z_2)\in\mathbb C^2\mid\mathfrak{Im}z_1>|z_2|^2\}$ equipped with its Bergman metric. For $\epsilon>0$ small, let us consider a subdomain $M_\epsilon\subset\Omega$ given by
\[M_\epsilon=\{(z_1,z_2)\mid y_1>x_2^2+y_2^2/ \epsilon^2\}\]
\noindent
on which the matrix group $G$ consisting of matrices
\[ g=\left[
\begin{array}{ccc}
\lambda^2 & 2i\lambda\xi & t+i\xi^2\\
    0   &\lambda  & \xi \\
    0   & 0       & 1       
\end{array}
\right]\qquad (\xi,\ t\in\mathbb R,\ \lambda>0)\]
\noindent
acts by matrix multiplication on $(z_1,z_2,1)^\mathrm{T}$. Easy computations show that for each $\epsilon>0$, $M_\epsilon$ is the strongly pseudoconvex total space of a $G$-bundle with $\dim_{\mathbb R}G=3$ acting in $M_\epsilon$ by holomorphic transformations and $\bar M/G\cong [-1,1]$. It turns out that $L^2\mathcal O$ on this space is trivial ({\it cf.} Section 3 of \cite{GHS}). 

\begin{rem}{\rm The examples in \cite{GHS} both have the property that the orbits of the singularity of the Levi function fill the boundary. It seems to us that amenability is impossible in this situation.}\end{rem}
\subsubsection{The thickened cylinder} Here we try to get amenability by taking the Cartesian product with a circle, {\it not} considered a group.
Consider the thickening of $(S^1\times\mathbb R)$,  
\[(S^1\times\mathbb R)^{\mathbb C}\cong \mathbb C/\mathbb Z \times \mathbb C\subset\mathbb C^2.\]
With $z_k = x_k+iy_k$, ($k=1,2$) the defining function $\rho:(S^1\times\mathbb R)^{\mathbb C}\to\mathbb R$ we can take to be $\rho(z_0,z_1) = y_0^2+y_1^2-1$. In particular, $\rho$ is invariant under translations 
\[\mathbb R\ni t:(z_0,z_1)\longmapsto(z_0,z_1+t).\]
Put $\mathcal T = \{(z_0,z_1)\in(S^1\times G)^{\mathbb C}\mid \rho(z)<1\}$. Then $\mathcal T\cong S^1_{x_0}\times\mathcal B^2_{(y_0,y_1)}\times\mathbb R_{x_1}$ and the quotient $\mathcal T/G\cong S^1\times\mathcal B^2$, the solid torus, 
\[\mathbb R\longrightarrow\mathcal T\longrightarrow S^1\times\mathcal B^2.\]
So pick a $\xi\in bM$ with coordinates $(\tilde x_0,\tilde y_0,\tilde x_1,\tilde y_1)$, ($\tilde y_0^2+\tilde y_1^2=1$) and consider its orbit: $G\cdot\xi=\{(\tilde x_0,\tilde y_0,t,\tilde y_1)\mid t\in\mathbb R\}\subset b\mathcal T$.
The Levi form is $L=2\delta_{jk}$, $(jk=0,1)$ and the Levi polynomial at the point $p\in bM$ with coordinates $\tilde\xi = (\tilde x_k +i\tilde y_k)_0^1$ is 
\[f(z,\tilde\xi)=-2i\sum_0^1\tilde y_k (z_k-\tilde\xi_k)-\sum_0^1 (z_k-\tilde\xi_k)^2.\]
Now choose $p\leftrightarrow\tilde\xi=(i,0)$, {\it i.e.} $y_0=1$, so that 
\[f(z)=-2i (z_0-i)-(z_0-i)^2-z_1^2.\]
\noindent
and the group only translates $z_1$. A point on the orbit of the singularity will be $p \leftrightarrow z=\tilde\xi$ so we take the convolution and take a limit as we approach that point:
\[(R_\Delta \chi f^\tau)(z) = \int_G dt\ \Delta(t)\chi(zt) \left[-2i (z_0-i)-(z_0-i)^2- (z_1+t)^2\right]^\tau\]
\noindent
Put $r=-2i (z_0-i)-(z_0-i)^2$ and note that as $z\to \tilde\xi$, $r\to 0$. In terms of $r$, the convolution is 
\[(R_\Delta \chi f^\tau)(z) = \int_G dt\ \Delta(t)\chi(zt) \left[r-(z_1+t)^2\right]^\tau.\]
\noindent
Coming toward the base point (and thus the path of the singularity) from inside the manifold, we can take $p_r$ with coordinates $z_1=0$ and $r\to 0$:
\[(R_\Delta \chi f^\tau)(p) =\lim_{r\to 0} \int_G dt\ \Delta(t)\chi(p_r t) \left[r- t^2\right]^\tau\approx\lim_{r\to 0} \int_{-\epsilon}^\epsilon dt\ \left[r- t^2\right]^\tau\]
\[\approx\lim_{r\to 0} \int_0^\epsilon dt\ \left[r- t^2\right]^\tau\]
\noindent
which will diverge for $\tau\le-1/2$. For $\tau>-1/2$ or $\tau =0$ (meaning the logarithm), then the integral converges. Thus $\mathcal T$ is $\tau$ amenable for $\tau\le-1/2$ and all points of the boundary are base points for singular functions. Note also that this action is also amenable. 
%
%%%%%%%%%%%%%%%%%%%
\subsection{Reduction of the structure group}
%%%%%%%%%%%%%%%%%%%
%
We have seen that in some cases it is possible to make amenable actions from those that are not by reducing the structure group. Sometimes it may be possible to construct holomorphic functions with stronger singularities at the boundary by a similar reduction. It is thus in our interest while solving the Levi problem on a $G$-bundle $M$ to choose the structure group $H\subset G$ with dimension as small as possible. To this end we note that if $H\subset G$ then $\square$ is $H$-invariant too, but not necessarily $H$-Fredholm unless $G/H$ is compact ($H$ is unimodular by Theorem 8.36 \cite{Kn}). In this case (as in the reduction in the Expression \eqref{reduce}) we may profit by working the problem in the form $H\to M\to X\times (G/H)$ instead of $G\to M\to X$. 
%
%%%%%%%%%%%%%%%%%%%
\section{Main Theorem}
%%%%%%%%%%%%%%%%%%%
%
Assume that the action of $G$ on $M$ is amenable and choose $f=f_\tau$, the Levi function at $x\in bM$ with $R_\Delta\chi f\notin C^\infty(\bar M)$. Consider the Fourier expansion of $h=\bar\partial\chi f$, $h=\sum_k H_k\otimes\psi_k$ with $(\psi_k)_k$ an orthonormal basis for $L^2(X)$ as in \S \ref{actions}. Also, choose an $l\in\mathbb N$ for which $R_\Delta H_l\neq 0$ and let $\mathfrak D_\delta=\mathfrak D_\delta^l$. Further, let
\[\langle\langle\chi f\rangle\rangle_\delta = \{R_\Delta \chi f\mid \Delta\in\mathfrak D_\delta\}.\]
\noindent
Since $\mathfrak D_\delta\subset L^2(G)$, Corollary \ref{L2} and Remark \ref{slice} imply that $\langle\langle\chi f\rangle\rangle_\delta\subset L^2(M)$. Furthermore, $\chi f\in L^2(M)$ is in the domain of the Hilbert space operator $\bar\partial_{\mathcal H}$ ({\it cf.} \cite{Per} or \cite{FK}). Amenability guarantees that $\langle\langle\chi f\rangle\rangle_\delta\cap C^\infty(\bar M)=\{0\}$. 
\begin{lemma}\label{lemma} The restricted antiholomorphic exterior derivative $\bar\partial:\langle\langle\chi f\rangle\rangle_\delta\to\langle\bar\partial\chi f\rangle_\delta$ is a bijection.\end{lemma} 
\begin{proof}Since $\bar\partial$ is elliptic with analytic coefficients, its kernel contains only analytic functions. The small $X$-support of the members of $\langle\langle\chi f\rangle\rangle_\delta$ imply therefore that kernel of $\bar\partial|_{\langle\langle\chi f\rangle\rangle_\delta}$ is trivial. Since $R_\Delta \chi f$ is smooth in the interior of $M$ for $\Delta\in\mathfrak D_\delta$, we have $\bar\partial R_\Delta \chi f=R_\Delta \bar\partial\chi f$. Since $\langle\langle\chi f\rangle\rangle_\delta\subset L^2(M)$, $\bar\partial$ and $\bar\partial_{\mathcal H}$ coincide there.\end{proof}
\begin{theorem}Assume that the action of $G$ in $M$ is amenable. Then the space $L^2\cap\mathcal O(M)$ is infinite-$G$-dimensional.\end{theorem}
\begin{rem}{\rm The method is similar to using a Friedrichs mollifier on the equation $\square u = \bar\partial\chi f$. The group convolution $R_\Delta \square u = R_\Delta \square u = \bar\partial R_\Delta\chi f$, $\Delta\in C^\infty(G)$, by invariance.}\end{rem}
\begin{proof} Let $f=f_\tau$ be a function with properties verifying amenability of the action. Theorem 6.6 of \cite{Per} provides that the operator $\square$ ({\it cf.} \eqref{KL}) on its domain is $G$-Fredholm. 
\noindent
Lemma \ref{bigspaces} allows us to conclude that, for $\delta>0$ sufficiently small, there exist closed, invariant subspaces 
\[L_\delta\subset{\rm im}\square\cap\langle\bar\partial\chi f\rangle_\delta\]
\noindent
of arbitrarily large $G$-dimension. For such $\delta$, let $0\neq g\in L_\delta$ and solve 
\[\square u = g.\]
\noindent
By Lemma \ref{smoothbig}, $\langle\bar\partial\chi f\rangle_\delta\subset C^\infty(\bar M,\Lambda^{0,1})$, so $g$ is smooth. The regularity of $\square$ then gives that $u\in C^\infty(\bar M)$. Since the images of $\bar\partial$ and $\bar\partial^*$ are orthogonal, we have 
\[\square u = \bar\partial\bar\partial^* u = g\] 
\noindent
and $g= \bar\partial\phi$ for some $\phi\in \langle\langle\chi f\rangle\rangle_\delta$ by Lemma \ref{lemma}.  Form the holomorphic function 
\[ \Phi = \phi - \bar\partial^* u.\]
\noindent
Amenability gives $\langle\langle\chi f\rangle\rangle_\delta\cap C^\infty(\bar M)=\{0\}$, from which $\phi\notin C^\infty(\bar M)$. We conclude that $\Phi\notin C^\infty(\bar M)$ and thus is nonzero.\end{proof} 
\noindent
The holomorphic function $\Phi$ in the proof cannot be extended smoothly beyond $bM$.
\begin{corollary} In the setting above, let $x\in bM$ be the base point of the Levi polynomial $f$. Then there exists a holomorphic function $\Phi_x$ which cannot be holomorphically extended beyond $x$.\end{corollary}
\begin{proof}Obvious.\end{proof}

%
%%%%%%%%%%%%%%%%%%%
\section{Acknowledgments} 
%%%%%%%%%%%%%%%%%%%
%
The author wishes to thank Gerald Folland, Jonathan Rosenberg, Matt Stenzel, and Alex Suciu for numerous conversations and Mikhail Shubin for the suggesting the problem and many years of friendship and mathematical advice.
\end{document}